%% file: pl.tex
\input amstex
\documentstyle {amsppt}
\tolerance=3000
\openup 6 pt
\nologo

\bigskip
\bigskip
\bigskip
\bigskip
\topmatter
\title
Inducing, Slopes, and Conjugacy Classes
\endtitle
\author
Roza Galeeva\footnote {Dept. of Mathematics, Northwestern University, Evanston, IL 60208-2730.}, Marco Martens \footnote {Institute of Mathematical Sciences, SUNY at Stony Brook, Stony Brook, NY 11794-3651.}, and Charles Tresser\footnote {I.B.M., Po Box 218, Yorktown Heights, NY 10598.}
\endauthor
\endtopmatter
\bigskip
\bigskip
\bigskip
\bigskip
\bigskip
\bigskip
\centerline {\bf Abstract.}
\bigskip
We show that the conjugacy class of an eventually 
expanding continuous  
piecewise affine interval map is contained in a smooth codimension 1 
submanifold of parameter space. In particular conjugacy classes have  empty interior. This is 
based on a study of the relation between induced Markov maps
and ergodic theoretical behavior.
\newpage

\bigskip
\centerline{\bf 1. Introduction}
\bigskip

\flushpar
One of the central questions in iteration theory is to decide whether two 
maps $f:X\to X$ and $g:Y\to Y$ are topologically conjugate, {\it i.e.}
 whether 
there exists a homeomorphism $h:X\to Y$ such that $h\circ f=g\circ h$.
In this paper we deal with this question  for multimodal 
continuous piecewise affine maps on the interval. A map is called
{\it piecewise affine} if it is continuous, piecewise monotone, and 
affine on each interval of monotonicity.

\flushpar
The unimodal case was 
already studied by Brodiscou, Gillot and Gillot [BGG]. They dealt with a 
one parameter family of unimodal piecewise
affine maps where just one slope varies. Later Misiurewicz and Visinescou
showed that the conjugacy classes in the general unimodal piecewise affine
family, where both slopes are allowed to vary, form lines in the parameter 
plane [MV].
In this paper we will show that the conjugacy class of an eventually
expanding piecewise affine map is contained in a codimension 1 submanifold
of parameter space.

\flushpar
The following trivial remark is the key to our study of the conjugacy
classes of piecewise affine maps. Consider a set $A$, finite or 
denumerable,
and assume that the interval $I$ is, up to a set of measure zero, the
pairwise disjoint union of intervals $I_a$, $a\in A$. Then a  map
$f:\bigcup I_a\to I$, where $f$ caries each $I_a$ onto $I$ by an affine 
homeomorphism,
is called a {\it multiple covering map (with index set $A$)}: notice that the domain of such a map has full measure in its image (in short, {\it has 
full measure}). The derivative of the branch $f|I_a$ 
is denoted by $Df_a$. 

\proclaim{Multiple Covering Map Principle} Let $f$ be a multiple-covering map
with index set $A$. Then
$$
\Sigma_{a\in A} \frac{1}{|Df_a|}=1.
$$ 
In particular, for multiple covering maps $f,\,g$ which have  
the same index set $A$:
$$\lbrace \forall a\in A,\,  
|Df_a|\ge |Dg_a|\rbrace\,\,\Rightarrow\,\,
\lbrace\forall a\in A,\,|Df_a|=|Dg_a|\rbrace\,.
$$
\endproclaim  

\flushpar
For piecewise affine Markov maps, the above principle applies almost
immediately, yielding similar results. As we shall see, a much finer 
analysis is required to deal with more general piecewise affine maps.
The first step of our study is to associate induced Markov maps (see \S 4) 
to piecewise 
affine maps. These induced maps have a topological 
definition and look like multiple covering maps. The only difference is that 
we don't know whether, for such a map, the domain of definition has full measure. So, before
applying the multiple covering map principle to induced Markov maps, we have to study the measure of their domain of definition. 

\flushpar
A piecewise affine map has the {\it Markov property} 
if it has an induced Markov map whose 
domain of definition has full measure, {\it i.e.}, 
it is a multiple covering maps (see \S4).
A closed set $A\subset N$ is called an {\it absorbing} set of
the  interval map $f:N\to N$ if
$$
\{x\in N| \omega(x)\subset A\}
$$  
has positive Lebesgue measure (
where $\omega(x)$ denotes the positive $f$-limit set of $x\in N$).
In [M] it was shown that
$S-$unimodal maps have the Markov property if and only if the map does not
have zero-dimensional absorbing sets. 
The main part of this paper is devoted to prove the same
result for piecewise affine maps. In Section 4 we prove

\proclaim{Theorem A} A piecewise affine map has the Markov property if
and only if it does not have zero-dimensional absorbing sets.
\endproclaim

\flushpar
There are three different properties which allow a map to have a 
zero-dimensional absorbing set. A map can have a periodic attractor.
Secondly it 
can be infinitely renormalizable (see \S2). In this case the topological
structure causes an absorbing  Cantor set. Furthermore  
a non-renormalizable map can have an absorbing Cantor set, which is caused
by intrinsic geometrical properties.

\flushpar
In [LM] and [L] it was
shown that quadratic unimodal maps can only have an absorbing  Cantor set if
they are infinitely renormalizable. However
recently it has been shown in [BKNS] that there exist unimodal maps with 
highly degenerate 
critical point having  absorbing Cantor sets. These results 
depend on a fine control of the geometry.
We will avoid such geometrical studies by an ergodic theoretical
shortcut: we only study eventually expanding piecewise affine maps, and
for such maps, one knows the existence of absolutely continuous
invariant probability measures [LY]. 
Yet, these measures cannot coexist with zero-dimensional absorbing 
sets and we get
 
\proclaim{Theorem B} Eventually expanding piecewise affine maps have the
Markov property.
\endproclaim

\flushpar
However the study of the intrinsic geometry is just postponed. To prove
the following conjecture one would probably have to go into geometrical
considerations.

\proclaim{Conjecture} A piecewise affine interval map with no periodic 
attractor is eventually expanding.
\endproclaim

\flushpar
In Section 5 we will apply the multiple covering map principle to
the Markov maps which are now, by Theorem B, multiple covering maps.  
A {\it branch} of a map is the restriction to an interval of monotonicity.
Branches of a piecewise affine map 
which contain pieces of the non-wandering set in the interior of its domain
are called {\it essential branches}.
 Let $\Cal{E}_d$ be the family of 
d-modal eventually expanding piecewise affine maps. $\Cal{E}_d$ is naturally
parametrized by some submanifold of ${\Bbb R}^{2d+2}$, to which we identify it.
Studying the conjugacy problem in $\Cal{E}_d$ we got the following result.

\proclaim{Theorem C} 
Every conjugacy class is contained in a codimension 1 submanifold of 
$\Cal{E}_d$.
Furthermore, if the slope of some essential branches of $f\in \Cal{E}_d$ 
are increased, the topological type changes.
\endproclaim

\flushpar
Section 2 contains some basic topological lemmas, some of which are part of the folklore. To simplify the exposition of the
proofs,  we only considered maps on the interval: most of this paper, and 
in particular Theorems A,B and C, hold true as well for piecewise affine circle maps which have at least one periodic orbit. In Section 3 we define {\it good} intervals and describe their
properties: they are the main ingredient in the definition of the Markov maps.

\flushpar
The proof of the characterization of maps with the Markov property,
presented in Section 4, 
can easily be generalized by using the tools from [M]. Thus Theorem A also
holds for smooth multimodal maps with 
negative Schwarzian derivative. More work would be needed to get a similar 
$C^2$ result.

\bigskip
\flushpar
{\bf Acknowledgement.} We thank John Milnor for suggesting improvements
to this text. He also raised the question as to whether 
isentrops, {\it i.e.}, parameter sets 
corresponding to maps with a given entropy, have empty interior.

\bigskip
\centerline{\bf Notations}
\bigskip

\flushpar
We will use the following conventions and notations. {\it Intervals} will always have positive length. Let $N=[a,b]$ be an interval 
and $A,U\subset N$,
where $A$ is measurable and $U$ open. Let $f:N\to N$ be a piecewise affine map.
\parindent=15pt
\item{-} $\partial U$ is the boundary of $U$,
\item{-} $int(A)$ is the interior of $A$,
\item{-} $mesh(U)$ is the length of the longest connected component of $U$,
\item{-} $|A|$ is the Lebesgue measure of $A$,
\item{-} $C_f$ is the set of critical points of $f$,
\item{-} $orb(A)=\{A,f(A),f^2(A),\dots\}$ is the orbit of $A$,
\item{-} $Df_i$ is the derivative of the $i^{th}$ branch.

\flushpar
We shall also say that  $U\subset N$ {\it satisfies the $\partial-$condition} 
if $orb(\partial U)\cap U=\emptyset$. A {\it branch } of a piecewise  
monotone map is the restriction
of the map to a maximal interval on which it is monotone.

\input pl2.tex

\input pl3.tex

\input pl4.tex

\input pl5.tex

\input plref.tex

\bye

%% file: pl2.tex
\tolerance=3000

\bigskip
\centerline{\bf 2. Non Renormalizable Maps}
\bigskip

Let us begin with some definitions.

\bigskip

The continuous map $f:N\to N$ is {\it piecewise affine} if there exist points
 $a=a_0<a_1<\dots<a_d<a_{d+1}=b$
such that $f_i=f|{[a_i,a_{i+1}]}$ is affine, and $Df_iDf_{i+1}<0$. The points
$a_1, a_2,\dots,a_d$ are called {\it critical points}. We say the map is 
{\it $d$-modal} when we want to stress the number of its critical points, {\it multimodal} when $d\geq 1$ and {\it unimodal} when $d=1$.

\flushpar
Consider a piecewise affine map $f:N\to N$. With $I\subset N$ an interval and $n\ge 1$, the pair $(I,n)$ is called a
 {\it renormalization} of $f$ if 
\parindent=15pt
\item{-} $f^n(I)\subset I$ and $\overline{I}\ne N$,
\item{-} the interiors of $f^i(I)$, $i=0,\dots, n-1$ are pairwise disjoint.

\flushpar
A map which has a renormalization is called {\it renormalizable}. The orbit
$\bigcup f^i(I)$ is called a {\it cycle (with period $n$)}.  
A cycle is called {\it minimal} if $f^n|I$ is non-renormalizable. A pair $(I,n)$ is called a {\it trap} of the piecewise affine map 
$f:N\to N$  if $f^n(I)\subset I$ and $\overline{I}\ne N$.

\flushpar
The following two properties of non-renormalizable maps will be used over and
over again.

\proclaim{Lemma 2.1} Let $f:N\to N$ be non-renormalizable 
and piecewise affine. 
Then 
\parindent=15pt
\item{1)} $\partial N\subset orb(C_f)$,
\item{2)} $f^{-1}(x)\setminus (\partial N\cup C_f)\ne \emptyset$ for every 
$x\notin \partial N$,
\item{3)} for every interval $I\subset N$
$$
\bigcup_{i\ge 0} f^i(I)=N.
$$
\endproclaim

\demo{Proof} The proof of $1)$ and $2)$ is easily supplied and we proceed with
the proof of $3)$.

\flushpar
Let $I\subset N$ be an interval. Observe that $f$ cannot have 
periodic attractors. In [MMS] it was proved that a piecewise affine map $f$ 
without periodic attractor cannot contract 
intervals too much: $\inf_{i\ge 0}|f^i(I)|>0$. 
This implies that the connected components of $\bigcup_{i\ge 0} f^i(I)$ 
have a definite size. Hence the set can have only a finite number of connected
components. These components are permuted by $f$.
In particular they are eventually periodic. This gives rise to
a renormalization. Hence there is exactly one component which is dense in 
$N$. Clearly this invariant component contains $C_f$. Hence it 
contains $orb(C_f)\supset \partial N$, so it is $N$.

\hfill\hfill\qed $\,\,$ (Lemma 2.1)
\enddemo

\proclaim{Lemma 2.2} A non-renormalizable piecewise affine map does not
have traps.
\endproclaim

\demo{Proof} Let $f:N\to N$ be a non-renormalizable piecewise affine map.
Observe that every non-renormalizable map has at least one expanding fixed 
point in $int(N)\setminus C_f$, say $f(p)=p$ and $|Df(p)|>1$.

\flushpar
Assume that there is some trap $(I,n)$. By Lemma 2.1 
$p\in I$. 

\flushpar
- Assume 
$p\in int(I)$ or $p\in \partial I$ is order preserving.
Let $E_k=\{x\in I|\{x,f(x),\dots,f^{k-1}(x)\subset I
\text{ and } f^k(x)\notin I\}$. 
Because there are no renormalizations $f(I)\not\subset I$. Hence $E_1\ne
\emptyset$. Assume $E_k\ne\emptyset$ for some $k\ge 1$.
We set $F_k=I\setminus (E_1\cup E_2\cup\dots\cup E_k)$. From the assumption on
$p$, $F$ contains an interval, hence by Lemma 2.1 $f(F_k)\not\subset F_k$.
By definition $f(F_k)\subset I$, thus $f^{-1}(E_k)\cap I=E_{k+1}\ne\emptyset$. 

\flushpar
Now consider $E_n$: $E_{n}\ne\emptyset$ and $f^{n}(E_{n})\cap I=\emptyset$. 
But $I$ is a trap, so $f^n(E_n)\subset f^n(I)\subset I$, a contradiction.

\flushpar
- Assume $p\in \partial I$ is an order reversing
fixed point. Then $p\notin int(f^i(I))$ for $i\ge 0$. Now $I\cup f(I)$ is also a trap
but containing $p$ in its interior and we are back to the previous case, which
implies $I\cup f(I)\supset int(N)$.
But now, $I$ and $f(I)$ are the components of $N-\{p\}$ and
$f$ interchanges these two components. We found a renormalization, a
contradiction.

\hfill\hfill\qed $\,\,$ (Lemma 2.2)
\enddemo

\proclaim{Lemma 2.3} The periodic points of a non-renormalizable 
piecewise affine map $f:N\to N$ are dense.
\endproclaim

\demo{Proof} Fix an open interval $I\subset N$ with 
$\overline{I}\cap \partial N=\emptyset$. The aim is to show that there exists 
$k\ge 1$ such that $f^k(I)\supset I$: $I$ contains a periodic point.

\flushpar
Because the orbit of $I$ is dense there exists a $q\ge 1$ such that
$f^q(I)\cap I\ne \emptyset$. Let $I_j=f^{jq}(I)$, $j\ge 0$. Assume that
$I\subset I_j$ never happens. 

\flushpar
Consider $T_k=\cup_{j\le k} I_j$, $k\ge 0$. Because $I_0\cap I_1\ne\emptyset$
every $T_k$ is an interval. Clearly $f^q(T_k)\subset T_{k+1}$ and 
$T_k\subset T_{k+1}$. Let $I=(a,b)$ and suppose $a\notin I_1$. 

\proclaim{Claim}
$$
a\not \in T_k \text{  for } k\ge 0.
$$
\endproclaim

\demo{Proof of Claim} For $k=0,1$ the claim is true. Assume by contradiction 
that there exists a first $k\ge 1$ such that $a\in T_{k+1}$.
Because $f$ is non-renormalizable, by Lemma 2.2  we get that 
$T_l-T_{l-1}=J_l\ne\emptyset$ 
for $l\le k$. Observe that $J_l\subset I_l$, hence $J_l$ is an interval,
otherwise $a\in T_l$. Now there exists $x\in \overline{J_{k-1}}$ 
with $f^q(x)=b_k$ where $b_k$ is the right boundary point of $T_k$. Let
$D=[x,b_{k-1}]$. Then $D\subset J_{k-1}\subset I_{k-1}$.

\flushpar
We have $a\in f^q(J_k)$ and $f^q(x)=b_k$, hence $T_k\subset f^q(D\cup J_k)$.
Now $I_0\not\subset f^q(J_k)$ by assumption, hence 
$$
D\cup J_k\subset T_k\setminus I_0\subset f^q(D).
$$
This yields 
$$
T_k\subset f^q(J_k\cup D)\subset f^{2q}(D),
$$
hence 
$$
I_0\subset f^{2q}(I_{k-1})=I_{k+1},
$$
a contradiction.

\hfill\hfill\qed $\,\,$ (Claim)
\enddemo

\flushpar
To finish the proof of Lemma 2.3, let $T=\cup T_k$. Now $f^q(T)\subset T$
and the closure of $T$ is not the whole $N$ because $T$ lies on one side of 
$I$ and $I$ does not touch the boundary of $N$. We found a trap, a 
contradiction. Hence for some $j\ge 1$ we have $I\subset I_j$.

\hfill\hfill\qed $\,\,$ (Lemma 2.3)
\enddemo

\proclaim{Corrolary 2.4} The critical set of a non-renormalizable piecewise
affine map has a neighborhood $U$ satisfying the $\partial-$condition,
$orb(\partial U)\cap U=\emptyset$,  
and
having arbitrarily small mesh.
\endproclaim
 
\proclaim{Lemma 2.5} Let $c\in C_f$ be a critical point of the 
non-renormalizable piecewise affine map $f:N\to N$. Then in every component 
$M\subset N\setminus (\partial N\cup C_f)$ there exists an open interval $I\subset M$ and
$n\ge 1$ such that $f^n|I$ is monotone and $c\in f^n(I)$.
\endproclaim

\demo{Proof} Fix $c\in C_f$. Let $d\in C_f$. If $c\in orb(d)$ then  
$n_d\ge 0$ will stand for the first time that $d$ hits $c$. 
Now take $L> \max\{n_d|n_d<\infty\}$.

\flushpar
Using Lemma 2.1 (2) we can choose a sequence $c_0=c,c_{-1},c_{-2},\dots$
such that $f(c_{-(i+1)})=c_{-i}$ and $c_{-i}\notin C_f\cup \partial N$.
Then consider $c_{-L}$. By Lemma 2.1 (3) there is some $n\ge 1$ with
$c_{-L}\in f^n(M)$. In particular there exists an interval $J_1\subset M$
such that $f^n|J_1$ is monotone and $c_{-L}\in f^n(J_1)$. 
Choose $J_1$ to be maximal, which means that $f^n|J_1$ is a branch,
{\it i.e. }, by Lemma 2.1 (1) we know that $\partial f^n(J_1)\subset orb(C_f)$.

\flushpar
Assume that $c_{-L}\in \partial f^n(J_1)$: then there would be some critical 
point
$d\in C_f$ and some $i\ge 1$ such that $f^i(d)=c_{-L}$. Take the pair
$(d,i)$ with $i$ minimal. Then $n_d=L+i>L$, contradicting the definition
of $L$. Hence $c_{-L}\in int(f^n(J_1))$. 

\flushpar
To finish the proof, consider the orbit of $c_{-L}$. It does not pass trough
critical points. Hence there is some open interval $J_2\ni c_{-L}$
with $f^L|J_2$ monotone and $c\in f^L(J_2)$. And we can take
$I=f^{-n}(J_2)\cap J_1$.

\hfill\hfill\qed $\,\,$ (Lemma 2.5)
\enddemo

\proclaim{Lemma 2.6} Let $f:N\to N$ be a non-renormalizable piecewise affine
map. 
Then for every interval 
$I$ there exist an interval $J\subset I$ and $n\ge 1$ such that 
$f^n|J$ is monotone and $f^n(J)$ is a connected 
component of $N\setminus C_f$.
\endproclaim

\demo{Proof}
Consider the interval $I\subset N$.

\proclaim{Claim} There exists $n\ge 0$
such that $f^n(I)$ contains a component of $N\setminus C_f$.
\endproclaim

\demo{Proof of Claim}
Because $f$ does not have wandering intervals and periodic attractors,
there exist infinitely many $n\ge 0$ with $f^n(I)\cap C_f\ne
\emptyset$. In particular there exist a critical point $c\in C_f$ and 
non-negative numbers
$n,q$ such that $f^n(I)\ni c$ and $f^{n+q}(I)\ni c$. Now consider the
interval $T=\bigcup_{j\ge 0} f^{jq}(f^{n+1}(I))$. 
If either $f^n(I)$ or $f^{n+q}(I)$ contains two consecutive critical points,
we are done. Otherwise, because $f$ does not have traps,
 $f^{n+1}(I)\subset f^{q}(f^{n+1}(I))$ which implies 
$$
f^{jq}(f^{n+1}(I))=\bigcup_{i\le j} f^{iq}(f^{n+1}(I)).
$$
But $f^q(T)\subset T$ and, since $f$ does not have traps, we get $\overline{T}
=N$. By Lemma 2.1 (3) we have in fact $T=N$. Hence there exists $j_0\ge 1$
such that $f^{j_0}(I)=\bigcup_{i\le j} f^{iq}(f^{n+1}(I))=N$.

\hfill\hfill\qed $\,\,$ (Claim)
\enddemo

\flushpar
We are going to prove that $I$ contains an
interval which is mapped after some time monotonically onto a component
of $N\setminus C_f$. Assume $f^j(I)$ does not cover a component of $N\setminus C_f$ for
$j<n$ but $f^n(I)$ does contain a component.
For every $k\le n$ there exits $J_k\subset I$ such that
\parindent=15pt
\item{-} $f^k|J_k$ is monotone,
\item{-} $f^k(J_k)=f^k(I)$.

\flushpar
To prove this let $J_0=I$ and assume that $J_k$ is defined for some $k<n$. 
If $f|f^k(I)$
is monotone then let $J_{k+1}=J_k$. If $f|f^k(I)$ is not monotone then
there exists a unique critical point $c\in C_f$ with $c\in f^k(J_k)=f^k(I)$.
Let $L,R\subset J_k$ be the intervals which are mapped onto the two
components of $f^k(J_k)\setminus\{c\}$. We may assume that 
$f^{k+1}(L)\subset f^{k+1}(R)$. Now $f^{k+1}(I)=f^{k+1}(R)$ and
$f^{k+1}|R$ is monotone. So let $J_{k+1}=R$.

\flushpar
To finish the proof of Lemma 2.6, we just choose $J=J_n$.

\hfill\hfill\qed $\,\,$ (Lemma 2.6)
\enddemo

\flushpar
A piecewise affine map $f$ is called {\it eventually expanding} if there is 
an integer $n\ge 1$ so that $|Df^n|>1$ whenever this derivative is defined.

\proclaim{Lemma 2.7} Every eventually expanding piecewise affine map is 
non-renormalizable or has finitely many minimal cycles.
\endproclaim

\demo{Proof} Every cycle $(I,n)$ consists of pairwise disjoint intervals.
This implies that the number of critical points of renormalizations is
uniformly bounded. Hence there is always a branch of $f^n|I$, say $f^n:J\to I$
monotone and $J\subset I$ whose   
size is a definite fraction of $I$. But $|Df^m|\to \infty$ so that
 $f^n$ could not map
this piece into $I$ for $n$ big. We conclude that the period of the
renormalizations
is bounded for eventually expanding maps. Since each cycle contains at least 
one critical point we conclude that there are only finitely many minimal 
cycles.

\flushpar
To find a minimal cycle take a renormalization with maximal period; $(I,n)$.
Then there will be a smallest  interval $J\subset I$ such that $(J,n)$ is
still a renormalization. The orbit of $J$ is a minimal cycle.

\hfill\hfill\qed $\,\,$ (Lemma 2.7) 
\enddemo

\demo{Remarks} 

\flushpar
1) the interiors of minimal cycles are pairwise disjoint,

\flushpar
2) almost every point enters after some time a minimal cycle.
In particular every minimal cycle equals the conservative part of some
ergodic component.

\flushpar
3) All statements 2.1--2.7 remain true if piecewise affine is replaced
by ``continuous and with no homterval''.
\enddemo

%% file: pl3.tex
\tolerance=3000

\bigskip
\centerline{\bf 3. Good Intervals}
\bigskip

\flushpar
Fix a non-renormalizable piecewise affine map $f:N\to N$.
An open set $U\supset C_f$ is called a {\it nice neighborhood} of $C_f$
if it satisfies the $\partial-$condition and every connected component
contains exactly one critical point. We set
$U=\bigcup_{c\in C_f} U_c$. Corollary
2.4 states that there are nice neighborhoods $U$ with $mesh(U)$
arbitrarily small.

\proclaim{Definition 3.1} Let $U_c\subset U\subset N$ be a component of
the nice neighborhood $U$ of $C_f$. An interval $T\subset N$ is called a
{\it good interval (of time $n\ge 0$) for $U_c$} if $f^n:T\to U_c$ is
monotone and onto.
\endproclaim

\flushpar
Because every component of a nice neighborhood of $C_f$ contains a
critical point, every good interval has a well defined time needed for
reaching the nice neighborhood. The $\partial-$condition
implies easily that two intersecting good intervals $T_1$ and $T_2$
corresponding to the same nice neighborhood are nested: if
$T_1\cap T_2\ne \emptyset$ then either
$T_1 \subset T_2$ or $T_2\subset T_1$.
The following Lemma states that the collection of good intervals is
big.

\proclaim{Lemma 3.2} Let $U\subset N$ be a
 nice neighborhood of $C_f$ with
$mesh(U)$ small enough. For every critical
 point $c\in C_f$ there exists, in
every interval $I\subset N$, a good interval
 $T\subset I$ for $U_c$, whose  time
is at least $1$.
\endproclaim

\demo{Proof} Fix $c\in C_f$. Lemma 2.5 says that in every component $M$
of $N\setminus C_f$ there exists an open interval $I_M(c)\subset M$
 and $n_M(c)\ge 1$
such that $c\in f^{n_M(c)}(I_M(c))$ and $f^{n_M(c)}|I_M(c)$ is monotone.
Choose such an interval in every component $M$ of $N\setminus C_f$.
Let $V_c=\bigcap f^{n_M(c)}(I_M(c))$.
By Lemma 2.6 we will find, in every
interval $I\subset N$, an interval
$J\subset I$ and $n\ge 1$ such that $f^n:J\to V_c$ is monotone and onto.
Now Lemma 3.2 holds if we take $U$ small enough such that
$U\subset \bigcap V_c$.

\hfill\hfill\qed $\,\,$ (Lemma 3.2)
\enddemo

\flushpar
Observe that we can describe topologically how small $U$ has
 to be to apply
Lemma 3.2. In Section 5 we will discuss conjugacy classes. For this
we prefer to deal with topologically defined objects.

\flushpar
To avoid the annoying fact that the branches
 can
be restricted by the boundary points of $N$, we
 assume that the map $f$ is part of an
{\it extension}, {\it i.e.} there is a piecewise affine
 map $g:[-1,1]\to [-1,1]$ such that
\parindent=10pt
\item{-} $N\subset [-1,1]$,
\item{-} $g|N=f$,
\item{-} $g(\{-1,1\})\subset \{-1,1\}$,
\item{-} every point in $(-1,1)$ enters $N$ after some time.

\flushpar
The next Lemma explains why nice neighborhoods are nice.

\proclaim{Lemma 3.3} Let $U\supset C_f$ be a nice neighborhood.
If $f^i(x)\notin U$ for $i<n$ but $f^n(x)\in U_c\subset U$, there
exists a good interval $T\ni x$ of time $n$ for $U_c$.
\endproclaim

\demo{Proof}
Let $U=\bigcup_{c\in C_f} U_c$ be a nice neighborhood of $C_f$.
Take $x\in N$ such that $f^i(x)\notin U$ for $i=0,\dots,n-1$ and
$f^n(x)\in U_c$. Suppose by contradiction that
 $f^n(T)$ does not cover $U_c$,
where $T\ni x$ is the maximal interval on which
 $f^n$ is monotone. We assumed
$f$ to be part of an extension. Hence the monotonicity
 is restricted by some
critical point: there exist $i\le n-1$ and a critical
 point $d\in C_f$ such
that $d\in\partial f^i(T)$  and $f^{n-i}((d,f^i(x))\subset U_c$.
 By definition of
$n$ we know that $f^i(x)\notin U_d$. Hence
$(d,f^i(x))\cap \partial U_d\ne \emptyset$ which implies that
$orb(\partial U_d)\cap U\ne \emptyset$, a contradiction.

\hfill\hfill\qed $\,\,$ (Lemma 3.3)
\enddemo

\proclaim{Lemma 3.4} Let $U$ be a nice neighborhood for $C_f$.
 There exists a
closed set $\Lambda_U$ with
Lebesgue measure zero such that every component of the complement of
$\Lambda_U$ is a good interval for $U$.
\endproclaim

\demo{Proof}
Let $U=\bigcup_{c\in C_f} U_c$ be a nice neighborhood of $C_f$
and $\Lambda_U=\{x\in N|orb(x)\cap U=\emptyset\}$. Choose
a component $S$ of the complement of $\Lambda_U$ and $x\in S$.
 Assume that
the orbit of $x$
enters $U$ for the first time in $n$
steps, say $f^n(x)\in U_c$. This means that there exists a good interval
$T\ni x$ of time $n$. Observe that all $U_c$ are good intervals for $U$.
Furthermore
because of the $\partial-$condition we know that good
intervals are nested: $f^i(T)\cap U=\emptyset$ for $i<n$. Because
$f^n(\partial T)=\partial U_c$ the $\partial-$condition implies that
$orb(\partial T)\cap U=\emptyset$: $\partial T\subset \Lambda_U$
 and $T=S$. Hence
every component of the complement of $\Lambda_U$ is a good interval.

\flushpar
The orbits of points in $\Lambda_U$ stay outside
 the neighborhood $U$ of the
critical points. The fact that the Lebesgue measure of such sets is zero
is shown in [Mi] and [M].

\hfill\hfill\qed $\,\,$ (Lemma 3.4)
\enddemo

\proclaim{Corollary 3.5} Let $U_n$, $n=1,2,\dots$ be nice neighborhoods
of $C(f)$ with $mesh(U_n)\to 0$.

\flushpar
If $X$ is a forward invariant set with positive Lebesgue measure then
for every $n\ge 1$  there exists a component $C_n$ of $U_n$ such that
$$
\lim_{n\to\infty}\frac{|X\cap C_n|}{|C_n|}=1.
$$
\endproclaim

\flushpar
An {\it ergodic component} of $f$ is a forward and
 backward invariant set with minimal
positive Lebesgue measure. Corollary 3.5 shows that there are at most
as much ergodic components of $f$ as there are critical points.

\flushpar
Given a neighborhood of a critical point it will in general not satisfy
the nice property of Lemma 3.3. The next Lemma shows
 how we can deal with this
problem.

\proclaim{Lemma 3.6} Let $V\subset \text{int}(N)$ be
 an interval containing
one critical point $c\in C_f$ and satisfying the $\partial-$condition.
Let $K=\{x\in C_f| orb(x)\cap V\ne \emptyset\}$. Then there exists a
neighborhood $U=\bigcup_{d\in K} U_d$ of $K$ with the
 following properties.
\parindent=15pt
\item{1)} $U_c=V$ and every component $U_d$ of $U$ contains only one
critical point $d\in C_f$. Furthermore $U$ satisfies
 the $\partial-$condition,
\item{2)} there is a function $l:(0,1)\to {\Bbb R}$
 with $l(y)\to 0$ if $y\to 0$ such that $mesh(U)\le l(|V|)$,
\item{3)} the set $K$ is partitioned, say $K=\bigcup_{j\le s}K_j$
 so that
\item{-} $K_0=\{c\}$,
\item{-} for every $d\in K_j$, $j>0$, there exist an interval
 $T_d\ni f(d)$, $n\ge 0$ and $e\in K_{j-1}$ such that
  $f^n:T_d\to U_e$ is monotone and onto
and $U_d$ is the connected component of $f^{-1}(T)$ containing $d$,
\item{4)} if $f^i(x)\notin U$ for $i<n$ but $f^n(x)\in U_d\subset U$,
there exists an interval $T\ni x$ such that $f^n:T\to U_d$
is monotone and onto.
\endproclaim

\demo{Proof} We construct $U$ by induction.
Assume we defined the objects:
\parindent=15pt
\item{i} disjoint sets $K_0=\{c\},
 K_1, K_2,\dots,K_s\subset K\subset C_f$,
\item{ii} neighborhoods $W_0=V\subset
 W_1\subset\dots \subset W_s$ of the form
$W_l=\bigcup_{i=0}^l\bigcup_{d\in K_i} U_d$ where $U_d$ is the connected
component of $W_l$ containing $d$. In particular $W_0=U_c=V$ and
$W_l \supset \bigcup_{i=0}^l K_i$,
\item{iii} numbers $t_1,t_2,\dots,t_l$
 such that:
\parindent=25pt
\item{-} every point $d\in K_i$, $i>0$,
         enters for the first time $W_{i-1}$ after $t_i$ steps, say
         $f^{t_i}(d)\in U_e$ with $e\in K_j, j<i$,  
\item{-} there exists an interval $T\ni v=f(d)$ such that $f^{t_i-1}$ maps $T$
         monotonically onto $U_e$ and $f^{-1}(T)=U_{d}$,
\parindent=15pt
\item{iv} $W_i$, $i\le s$ satisfies the $\partial-$condition.

\flushpar
Assume that $\bigcup_{l\le s} K_l\ne K$. We are going to define
$K_{s+1}$, $W_{s+1}$ and $t_{s+1}$ according to the above properties.
Take $x\in K\setminus\bigcup_{i\le s} K_i$ and
 let $t_{s+1}(x)$ be the first moment
that the orbit of $x$ enters $W_s$. This happens because $x\in K$ and
$V\subset W_s$. Now let $t_{s+1}=\min\{t_{s+1}(x)\}$ and
$K_{s+1}=\{x\in K-\bigcup_{i\le s} K_i|t_{s+1}(x)=t_{s+1}\}$.
To finish the construction we have to find the intervals $U_d=f^{-1}(T_d)$
for the points $d\in K_{s+1}$. Choose $d\in K_{s+1}$, say
 $f^{t_{s+1}}(d)\in
U_e$ with $e\in K_j$, $j\le s$. Consider the maximal
 interval $M\ni v=f(d)$
on which $f^{t_{s+1}-1}$ is monotone and assume that
 the monotone image does
not cover $U_e$. Because we assumed that $f$ is part of an extension, the
monotonicity is restricted by some critical point and not
 by a boundary point.
There exits a critical point $e'\in C_f$ and a number
 $k<t_{s+1}-1$ such that
$e'\in \partial f^k(M)$ and $f^{t_{s+1}-1-k}((e',f^k(v)))$ is strictly
contained in $U_e$.
Observe that every  point in $W_s$ eventually enters $V$. So $orb(e')$
intersects $V$: $e'\in K$. Because $f^{t_{s+1}-1-k}(e')\in U_e$ and
$t_{s+1}-1-k<t_{s+1}$ we get $e'\in K_0\cup K_1\cup\dots\cup K_s$. Hence
$e'\in U_{e'}\subset W_s$. Because $k<t_{s+1}-1$ we have
$f^k(v)\notin U_{e'}$:
$\partial U_{e,}\cap (e',f^k(v))\ne \emptyset$. This implies that
$orb(\partial U_{e'})\cap U_d\ne \emptyset$. This cannot be because $W_s$
satisfies the $\partial-$condition. This contradiction implies
 $f^{t_{s+1}-1}(M)\supset U_e$.
Now we can take the interval $T_d\ni v$ which is mapped by $f^{t_{s+1}-1}$
monotonically onto
$U_e$, and we let $U_d$ be the connected component of $f^{-1}(T_d)$
which contains $d$. This finishes the definition of $W_{s+1}$.

\flushpar
To finish the induction step we have  to check that
 $W_{s+1}$ satisfies the
$\partial-$condition. To do so, take $y\in \partial
 W_{s+1}$ and assume by
contradiction that for some $n\ge 0$ we have $f^n(y)\in W_{s+1}$, say
$f^n(y)\in U_e$ with $e\in \bigcup_{j=0}^{s+1} K_j$.
Because every point in $W_{s+1}$ enters after
 some time $W_s$, we know that
$y\in \partial U_d$ with $d\in K_{s+1}$.
Because $f^{t_{s+1}}(y)\in\partial W_s$ and $W_s$ satisfies the
$\partial-$condition, we have $n< t_{s+1}$. Hence
$f^n(d)\notin U_e$. So $\partial U_e\cap f^n(U_d)\ne\emptyset$ and
$orb(\partial U_e)$ intersects $W_s$, a contradiction.

\flushpar		
This procedure will stop after finitely many steps:
$\bigcup_{j\le s} K_j=K$. Let $U=W_s$. Clearly $U$ satisfies the
$\partial-$condition.
The Contraction Principle from [MMS] implies that
$mesh(U)$ goes to zero if $|V|$ goes to zero.

\flushpar
It remains to prove property $4$. 
Take $x\in N$ and suppose that $x$ enters $U$
 for the first time in $n\ge 0$
steps, say $f^n(x)\in U_d$. Let $M\ni x$ be the maximal interval on which
$f^n$ is monotone and suppose that $f^n(M)$ does
 not cover $U_d$. Because we
assumed $f$ to be part of an extension, the monotonicity
 is restricted by a critical point: there exist $e\in C_f$
 and $i<n$ such that $e\in \partial f^i(M)$
and $f^{n-i}((e,f^i(x))$ is strictly contained in $U_d$.
 First observe that this implies $e\in K$.
So $f^i(x)\notin U_e$ which implies that $\partial U_e\cap (e,f^i(x))\ne
\emptyset$.
Hence $orb(\partial U_e)\cap U\ne \emptyset$. This is impossible because
$U$ satisfies the $\partial-$condition.

\hfill\hfill\qed $\,\,$ (Lemma 3.6)
\enddemo

%% file: pl4.tex
\tolerance=3000

\bigskip
\centerline{\bf 4. The Markov-Property }
\bigskip

\flushpar
As we will see in this section, every nice neighborhood of the critical
points of a piecewise affine map defines an induced map.
This induced map is strongly related to the ergodic theoretical behavior
of the map. In particular the existence of
absorbing Cantor sets is related
to these induced maps. We will start to define these induced maps for the
non-renormalizable piecewise affine  map $f:N\to N$.

\flushpar
Fix a nice neighborhood $U\subset N$
 of $C_f$ with $mesh(U)$ small enough to apply Lemma 3.2.
Let $D\subset N$ be the union
of all good intervals
for $U$ whose time is at least $1$. From Lemma 3.2 we know that $D$ is
dense in $N$, and because good intervals are disjoint
 or nested, we get that
every connected component of $D$ is a good interval of time at least $1$.
This allows us to define the {\it Markov map $M:D\to U$ (defined by $U$)}
in the following way: for every connected component
$T\subset D$ of $D$, with time $n$, set
$$
M|T=f^n|T.
$$
Observe that Markov maps are defined topologically.

\proclaim{Definition 4.1} A piecewise affine map $f:N\to N$ has
the {\it Markov property} if there exists a nice
 neighborhood $U\subset N$
of the critical points such that its
Markov map $M:D\to U$ is defined almost everywhere,{\it i.e.},
$$
|N-D|=0.
$$
The points in the set $B_0=N-D$ are called {\it bad points}.
\endproclaim

\flushpar
A closed set $A\subset N$ is called
 an {\it absorbing set} if
$$
|\{x\in N| \omega(x)\subset A\}|>0,
$$
where $\omega(x)$ denotes the $\omega$-limit set of $x\in N$. 

\proclaim{Theorem A} A non-renormalizable piecewise affine map has the
Markov property if and only if it does not have zero-dimensional 
absorbing sets.
\endproclaim

\flushpar
The next two lemmas are needed as preparation for the proof of Theorem A.
The first one gives a description of the limit behavior of bad points.
The second one is technical but will be used to prove
 the ergodicity of non-renormalizable maps.
 It also enables us to define special Markov maps whose
image is just one interval. These Markov maps play a crucial role
in the description of conjugacy classes in section 5.

\flushpar
To describe the limit behavior of bad points, we need some preparation.
Fix a nice neighborhood $U\subset N$ of $C_f$ and consider the Markov map
$M:D\to U$. Let $B_0=N-D$ be the closed zero-dimensional set of bad
points. In general
this set will not be invariant. The main step in the proof of Theorem A
is to find an ``almost invariant" set $\hat{B}\supset B_0$.

\flushpar
We start to define the {\it depth} $d(T)$ of a good interval $T$
as the number of good intervals which contain strictly $T$. Clearly every
component of $D$ has depth equal to zero.
If the critical value $f(c)$ is contained in a good interval of
 depth $d$,
denote this interval by $T_d(c)$. A critical point which has infinitely
many $T_d$'s is said to be of {\it infinite type}. Otherwise, it is
 of {\it finite type}.

\flushpar
We also need to define the {\it pull back} of bad points along the orbit
 of a good interval. Let $T\subset N$ be a good interval,
 say $f^n:T\to U_c$. Define the {\it tubes}
$$
P_T=\bigcup_{i=0}^n \overline{f^{-i}(B_0\cap U_c)\cap f^{n-i}(T)}.
$$
Now the {\it extended set of bad points} is defined to be
$$
\hat{B}=\bigcup_{n\ge 0} B_n,
$$
where
$$
B_n=B_0\cup\bigcup_{d=0}^n
           \bigcup_{c\in C_f} P_{T_d(c)}.
$$
Clearly every $B_n$ is a closed zero-dimensional set.

\proclaim{Lemma 4.2} The extended set of bad points contains
 $B_0$ and satisfies $f(\hat{B}-C_f)\subset \hat{B}$.
 For almost every $x\in B_0$, there exists $n\ge 0$ with
$$
orb(x)\subset B_n.
$$
\endproclaim

\demo{proof} Assume that $f$ is part of an extension.
First we will show the near invariance property of
 $\hat{B}$. Because the tubes $P_T$
are invariant, it suffices to show that $f(B_0-C_f)\subset \hat{B}$. Let
$x\in B_0-C_f$. If $f(x)\notin B_0$, there exist $c\in C_f$ and
$d\ge 1$ with $f(x)\in T_d(c)$. Because $x\ne c$, we can take $d$ maximal
with these properties. Now  $f(x)\notin \hat{B}$ implies
$f(x)\notin P_{T_d(c)}$. So there exists some good interval
 $T\subset T_d(c)$
with $f(x)\in T$, and because $d$ was taken to be maximal, we have
 $f(c)\notin T$.
Then $f^{-1}(T)$ contains a good interval around
 $x\in B_0$, a contradiction.

\flushpar
Now take some $x\in B_0$ and assume that $orb(x)\not\subset B_n$
for all $n\ge 1$. This means that after some time, $x$ hits
$V_d(c)=f^{-1}(T_d(c))\cap U_c$ for some $c\in C_f$
and some $d\ge n+1$.

\flushpar
There are two observations to be made:
\parindent=15pt
\item{-} $|V_d(c)|\to 0$ if $d\to \infty$,
\item{-} $T_d(c)$ never passes trough $V_d(c)$ (otherwise there would be an 
attractor).  This implies that $V_d(c)$
satisfies the $\partial-$condition.

\flushpar
This means that we can apply Lemma 3.6 and get nice extensions $W_d(c)$
of $V_d(c)$ with $mesh(W_d(c))\to 0$ if $d\to \infty$.

\flushpar
It only remains to show is that for some $\epsilon>0$,
$$
\frac{|D\cap W|}{|W|}\ge\epsilon
$$
for every component $W\subset W_d(c)$ and every $c\in C_f$ and $d\ge 1$:
we can push back this definite amount of good intervals 
into a very small neighborhood of $x$ by using Lemma 3.6 again, 
showing that $x$ is not a density point of $B_0$. 
Density points could not go to deep in $\hat{B}$, and the Lemma 
will be proved.

\flushpar
Let $K\subset \overline{K}\subset I$ be two open intervals. The set
$A=I-K$ is called a {\it boundary piece} if
$|\{T\subset A| T \text{ is a good interval }\}|=|A|.$
The first step is to show that every $U_c$ has a boundary piece. Fix
$c\in C_f$ and consider the sequence of intervals
$Q_1=f(U_c), Q_2,\dots$ with the following
 inductive definition: If $Q_i\subset T$
where $T$ is a good interval for $U_d$,
 then $Q_{i+1}=f(U_d)$. Otherwise the
sequence stops. If this sequence is longer
 than the number of critical points,
then at least one critical point is visited at least twice, and there is
 a trap.
Hence this sequence is finite. Say $Q_s=f(U_d)$ is not
 subset of $T_1(d)$.
Now apply Lemma 3.4 and we see that $U_d$ has a boundary piece.
 Considering
the sequence $Q_s,\dots,Q_2,Q_1$, we can pull back parts of this piece
and we will find a boundary piece in $U_c$.

\flushpar
The second step is to make definite boundary pieces in the $V_n(d)$.
 This is
easy because we can pull back one of the above boundary pieces into
$T_n(c)$, giving rise to definite boundary pieces in $T_n(d)$. One step
more and we will find the definite boundary pieces in  $V_n(c)$.

\flushpar
Lemma 3.6 describes how the different components of $W_n(d)$ are
related: they form a tree. Using this description we can pull back
the definite boundary pieces in $V_n(c)$ into definite boundary pieces of
the components of $W_n(c)$.

\flushpar
Observe that the only non-bounded  part of the construction takes place
during the transport of the boundary pieces in $U_c$ to the $T_n(d)$.
This transport is affine so that the proportion of space occupied
by boundary pieces is preserved, as well as the fact that
these boundary pieces are filled by good intervals.

\hfill\hfill\qed $\,\,$ (Lemma 4.2)
\enddemo

\proclaim{Lemma 4.3} Let $U=\bigcup_{c\in C_f} U_c\supset C_f$
be nice neighborhood with $mesh(U)$
small enough. The set $D_{\infty}$ consists of all points $x\in N$
which are contained in infinitely may good intervals:
$$
x\in  \dots \subset T_3(x)\subset T_2(x)\subset T_1(x)
$$
with $f^{t_i(x)}:T_i(x)\to U_{c_i(x)}$ and $t_i(x)\to\infty$.

\flushpar
For every critical point $c\in C_f$ and for almost all $x\in D_{\infty}$
there are infinitely $T_i(x)$ with $f^{t_i(x)}:T_i(x)\to U_c$.

\flushpar
In particular if $B\subset U_c$ with $|B_0\cap U_c-B|=0$ and $|B|>0$,
then almost every $x\in D_\infty$ hits $B$ after some time.
\endproclaim

\demo{Proof} Fix $c\in C_f$. Lemma 3.2 implies that every $U_d$ contains a
good interval for $U_c$. Let $B\subset U$ be the union of those good
intervals and $X_B=\{x\in D_\infty|orb(x)\cap B=\emptyset\}$. Take
$x\in D_\infty$ and consider the sequence
$x\in\dots\subset T_2(x)\subset T_1(x)$ of good intervals with times
$t_i(x)\to \infty$. Let $B_i=T_i\cap f^{-t_i}(B)$. Clearly
$B_i\cap X_B=\emptyset$. Because $f^{t_i}|T_i$ is affine we get
$$
\frac{|B_i|}{|T_i|}=\frac{|B\cap U_{c_i(x)}|}{|U_{c_i(x)}|}
\ge\min_{c\in C_f}\frac{|B\cap U_c|}{|U_c|}\ge \epsilon>0
$$
for all $i\ge 1$.
 Because there are no wandering intervals and no periodic
attractors we have $|T_i(x)|\to 0$.
 Hence $x$ is not a density point of $X_B$
and $|X_B|=0$. Now $D_\infty$ is covered up to a set of measure zero by
good intervals for $U_c$.
 From this we get directly that almost
 every point in $D_\infty$ is contained
 in infinitely many good intervals for $U_c$.

\flushpar
Take a set $B\subset U_c$ with $|B|>0$ and $|B_0\cap U_c-B|=0$. Let
$X_B=\{x\in D_\infty|orb(x)\cap B=\emptyset\}$. As above we can show that
$|X_B|=0$.
\hfill\hfill\qed $\,\,$ (Lemma 4.3)
\enddemo

\flushpar
In stead of proving Theorem A, we will prove the following stronger 
proposition describing more precisely the ergodic theoretical behavior. A map
$f:N\to N$ is called {\it ergodic}
 if it does not have two disjoint ergodic
components. It is called {\it conservative}
 if almost every point hits after
some time an arbitrarily given set $X\subset N$
 with positive Lebesgue measure.

\proclaim{Proposition 4.4} Let $f$ be a non-renormalizable
 piecewise affine
map.

\flushpar
If $f$ has the Markov property, $f$ is ergodic and conservative. In
particular the orbit of almost every point is dense in $N$.

\flushpar
If $f$ has a Markov map
$M:D\to U$ with $|N-D|>0$ then
there exists $s\ge 0$ such that for almost all $x\in N$
$$orb(f^{n_x}(x))\subset B_s
$$
for $n_x\ge 0$ big enough. In particular
$$
|N-\{x\in N|\omega(x)\subset B_s\}|=0,
$$
so that $B_s$ is a zero-dimensional absorbing set, 
absorbing in fact almost every orbit.
\endproclaim

\demo{Proof}
Let $f$ be a piecewise affine map
having the Markov property. The Markov property implies
 that $|N-D_\infty|=0$.
Now let $X\subset N$ be an invariant set of positive
 Lebesgue measure. Take
$c\in C_f$ and a density point $x\in X\cap D_\infty$ of $X$. Now consider
only the intervals
$T_i(x)$ from Lemma 4.3 which are good for $U_c$ and observe that
$\frac{|X\cap T_i|}{|T_i|}\to 1$.
Because $X$ is invariant we conclude that $\frac{|X\cap U_c|}{|U_c|}=1$.
Conclusion: we cannot have two disjoint invariant sets
 of positive measure: the map $f$ is ergodic.

\flushpar
To prove the conservativity of $f$, take a set
 $A\subset N$ with positive
Lebesgue measure. From Proposition 2.1 we know that there is some
$J\subset U_c$ with positive Lebesgue measure
and some number $n\ge 0$ such that $f^n(J)\subset A$.
Now apply Lemma 4.3 to $B=J$.
 Almost all point enters $J$ after some time,
hence also enters $A$ a little bit later.

\flushpar
Consider next a piecewise affine map $f$ which
does not have the Markov property. Then there exists a Markov map 
$M:D\to U$ with $|B_0|>0$. From Lemma 3.4
 we get $|B_0-U|=0$. Hence there exists a $c\in C_f$
with $|U_c\cap B_0|>0$. As a direct consequence of Lemma 4.3, we get
$|D_\infty|=0$: almost every point hits $B_0$ after some time.
The limit behavior of orbits is guided by the behavior of points
in $B_0$. This behavior is described by Lemma 4.2, giving rise to the
following candidates for the ergodic components
$$
E_n=\bigcup_{i\in\Cal Z} f^i(E'_n),
$$
where
$$
E'_n=\{x\in B_0|orb(x)\subset B_{n} \text{ and } orb(x)\not\subset
B_{n-1}\}.
$$
These sets $E_n$ are pairwise disjoint,
 forward and backward invariant sets.
Furthermore by Lemma 4.2, we get $|\bigcup_{n\ge 0}E'_n|=|B_0|$, and as a
consequence of Lemma 4.3, $|N-\bigcup_{n\ge 1} E_n|=0$. Now Corollary 3.5
implies that  that there are only finitely many $E_n$ with $|E_n|>0$.
Hence for some $s\ge 0$
$$
|N-\bigcup_{n\le s} E_n|=0.
$$
This means that the limit behavior takes place in $B_s$:
$$
|N-\{x\in N| \omega(x)\subset B_s\}|=0,
$$
so that $B_s$ is a zero-dimensional set absorbing almost all points in $N$.

\hfill\hfill\qed $\,\,$ (Proposition 4.4)
\enddemo

\demo{Remark} Proposition 4.4 implies that if a map has
 the Markov property,
then all its Markov maps are defined almost everywhere.
\enddemo

\proclaim{Theorem B} An eventually expanding non renormalizable piecewise
affine map has the Markov-Property.
\endproclaim

\demo{Proof} In [LY] it was proved that an eventually expanding
 piecewise affine map has an absolutely continuous invariant
 probability measure, and that furthermore,
 the density of this measure has bounded variation.

\flushpar
Now assume that there is an eventually expanding non renormalizable
 piecewise affine map not having the Markov property. Given a Markov map
$M:D\to U$, there is some $s\ge 0$ such
that the orbit of almost every point enters the closed set $B_s$
 after some time.
 Hence every ergodic component of the invariant probability
 measure is supported on
 $B_s$. In fact the whole measure is supported on $B_s$. Now $B_s$ is
 zero dimensional and closed,
 and such sets cannot support a
 non-zero density of bounded variation. This yields a
 contradiction.

\hfill\hfill\qed $\,\,$ (Theorem B)
\enddemo

%% file: pl5.tex
\tolerance=3000

\bigskip
\centerline{\bf 5. Conjugacy Classes}
\bigskip

\flushpar
In this section we are going to consider families of piecewise
 affine maps
and show that every conjugacy class is contained in a submanifold of
codimension 1 in the space of such maps.

\flushpar
Let $\Cal{F}_d$ be the family of $d-$modal piecewise affine maps. 
 The subfamily of
eventually expanding piecewise affine maps is denoted by
$\Cal{E}_d\subset \Cal{F}_d$. We consider $\Cal{E}_d$ as a submanifold 
of ${\Bbb R}^{2d+2}$.
We study the conjugacy question inside the class $\Cal{E}_d$. The conjugacy 
class of a map $f\in \Cal{E}_d$ is denoted by $[f]\subset \Cal{E}_d$.

\flushpar
To describe the conjugacy classes we need the notion of
essential branches and slopes. The $i^{th}$ branch is called
 {\it essential} if a minimal cycle intersects the interior of its domain.
 The slope
 $Df(i)$ is then also  called {\it essential}.
 Let $B_f$ be the collection of the
essential branches. Observe that $B_f$ is defined topologically.
In general one can change non essential slopes of
 a map without changing its topological type: examples are easily
provided.

\proclaim{Theorem C} The conjugacy class $[f]\subset \Cal{E}_d$ is
contained in a codimension 1 submanifold of $\Cal{E}_d$.
In particular if $g\in [f]$
and its essential slopes are at least as big as the corresponding
essential slopes of $f$ then they are in fact equal, {\it i.e.},
$$\{|Dg(i)|\ge |Df(i)|,\, i\in B_f\}\,\Rightarrow\,
 \{Df(i)=Dg(i)\,\,{\text for}\,\, i\in B_f\}.$$
\endproclaim

\flushpar
The proof of this Theorem is based on the
 Multiple Covering Map principle.
We will not work in $\Cal{E}_d$, but in the space of
 inverses of  slopes.
To go back to $\Cal{E}_d$ we use the submersion
$$
\pi:\Cal{E}_d\to D_d=(0,\infty)^{d+1}
$$
defined by $\pi(f)(i)=(|Df(i)|)^{-1}$.

\flushpar
The basic step in proving Theorem C is the definition of an induced map.
This will allow us to define topologically a multiple covering map  for 
every $f\in \Cal{E}_d$. 

\flushpar
Choose a map $f\in \Cal{E}_d$ and consider the minimal
 cycle corresponding
to a non renormalizable renormalization $(N,n)$.
Let $B\subset B_f$ be the collection of
essential branches which are used by the cycle $\bigcup f^i(N)$ and let
$D=(0,\infty)^B$. The natural projection from $D_d$ into $D$ is denoted
by $p:D_d\to D$.
Let $g=f^n|N$ and take a nice neighborhood $U\subset N$ of $C_g$. Choose
$c\in C_g$ and consider the union $G\subset U_c$ of all good intervals
$T\subset U_c$ for $U_c$ of positive time. As in the definition of
Markov maps we get an induced map $T:G \to U_c$.
Since we started with an eventually expanding map $f$, the map $g$ is also
eventually expanding. Hence it has the Markov property. Now by applying
Lemma 4.3 we get that
$$
|U_c-G|=0,
$$
{\it i.e.}, $T$ is a multiple covering map.

\flushpar
Before applying the Multiple Covering Map principle,
 we need some definitions.
Let $B_T$ be the collection of branches
of $T$. For every $I\in B_T$ there is a unique $t_I\ge 1$ such that
$T|I=f^{t_I}|I$. The only thing left over is to count how many times the
orbits of those branches of $T$ use the branches of $f$.
 Let $I\in B_T$ and
$i\in B$ and define $t_I(i)=\#\{j\le t_I-1|f^j(I)\subset i\}$.

\flushpar
Let $p\circ \pi(f)=(y_1,\dots, y_b)\in D$. Then the
Multiple Covering Map principle tells us
$$
\sum_{I\in B_T}\prod_{i\in B} y_i^{t_I(i)}=1.
$$
Now observe that the objects $U_c, T, B, B_T$ and $t_I(i)$ are all
topologically defined. So if we define a real analytic function
$\psi: D\to{\Bbb R}$ by
$$
\psi(x)=\sum_{I\in B_T}\prod_{i\in B} x_i^{t_I(i)}
$$
then $\psi\circ p\circ \pi(f')=1$ for all $f'\in [f]$. So
$$
p\circ \pi([f])\subset \psi^{-1}(1).
$$

\def\mar#1{\buildrel{\scriptstyle #1}\over \longrightarrow}
\flushpar
The sequence
$$
\Cal{E}_d\mar{\pi}D\mar{p}D\mar{\psi}{\Bbb R}\ni 1
$$
indicates how to prove Theorem C: we have to show that $1$ is a regular
value of $\psi\circ p\circ \pi$. In particular we will show that
the gradient of $\psi$ has only positive entries. This will also imply
the second statement of Theorem C.

\demo{Proof of Theorem C}
Take $f\in \Cal{E}_d$ and let $F=p\circ \pi([f])\subset D$.
Let $R>0$ be the radius of convergence of $\psi$. Then for $x\in D$
with $|x|<R$ the gradient of $\psi$ is defined and
$$
\frac{\partial\psi}{\partial x_i}(x)=
\sum_{I\in B_T(i)}t_I(i)x_i^{t_I(i)-1}
\prod_{j\ne i} x_j^{t_I(j)}>0,
$$
where $B_T(i)=\{I\in B_T| t_I(i)>0\}$, {\it i.e.} consists of
 those branches whose orbits
pass trough the branch $i\in B$. Such branches exist because $g|N$ is
non-renormalizable: the $g-$orbit of $U_c$ is dense in $N$,
 hence the $f-$orbit of $U_c$ is dense in the minimal cycle.

\flushpar
The fact that all components of the gradient of
 $\psi$ are positive implies
that $1$ is a regular value of $\psi$.
 So $\psi^{-1}(1)\cap\{x|\,|x|<R\}$
is an analytic codimension 1 submanifold of $D$.

\demo{case 1 - $F\subset \{x|\,|x|<R\}$} Then $F$
is part of the analytic codimension 1 submanifold
$\psi^{-1}(1)\cap \{x|\,|x|<R\}$, hence $[f]$ is
 contained in the analytic
codimension 1 submanifold $(\psi\circ p\circ \pi)^{-1}(1)$.
Furthermore, because $\psi$ has a positive
gradient, we will leave $\psi^{-1}(1)$, and so $F$,  by increasing some
essential slopes. Theorem C is proved.
\enddemo

\demo{case 2 - $F\not\subset \{x|\,|x|<R\}$ }
Because $\psi=1$ on $F$ we have $F\subset \{x|\,|x|\le R\}$.
 We will use polar
coordinates on $\{x|\,|x|\le R\}$. Let $S=\{x\in D|\, |x|=R\}$. Then
$\{x|\,|x|\le R\}=S\times (0,R]$.

\flushpar
Consider $L=\psi^{-1}(1)\cap (S\times (0,R])$.
Because the gradient of $\psi$ has only positive entries we
have $\#\{L\cap (\{\theta\}\times (0,R])\}\le 1$.
 The set $X\subset S$ consists
of those $\theta$'s whose rays intersect $L$.
Using again the positivity of the gradient of
$\psi$, we get that the projection $L\to X\subset S$ is a diffeomorphism.
The invariance of domain gives us that $X$ is open.

\flushpar
Now we can consider $L\subset S\times (0,R]$ to be the graph of
an analytic  function $\phi:X\to (0,R]$. Observe that $\phi(x_n)\to R$
whenever $x_n\in X$ converges to $x\in \partial X$. This follows from
the fact that $L$, as a preimage under $\psi|S\times (0,R)$, is closed in
$S\times (0,R)$ and the fact that $\partial X\cap X=\emptyset$.
This allows us to extend $\phi$ continuously to a map $\phi:S\to (0,R]$
by defining $\phi|S-X=R$.

\flushpar
Now $F\subset graph(\phi)\subset S\times (0,R]$. The $graph(\phi)$ is
a continuous codimension 1 submanifold of $D$. It remains to describe what
happens if we increase an essential slope. Clearly it is sufficient
 to show
that a decrease of some essential slopes will change the topological
 type. If our map $f$ is projected into $L$, then by decreasing
essential slopes we will
increase $\psi$ and hence leave $p\circ \pi([
f])$. On the other hand if
$p\circ \pi(f)=(\theta,R)\in S$, a decrease of some essential slopes
will cause that we leave $S\times (0,R]$ and so $\psi^{-1}(1)\subset
S\times (0,R]$.

\hfill\hfill\qed  $\,\,$ (Theorem C)
\enddemo
\enddemo

\flushpar
Markov maps can be used to define Markov extensions (the original map
is a factor of the extension). Those extension turned out to be very
useful for studying the ergodic theoretical behavior, specially 
the existence of invariant measures.
 This is because extensions
always have an absolutely continuous invariant measure.
 Invariant measures
for the original map could be constructed by trying to
 project the measure of the extension.
However the question whether every continuous measure can be obtained
by projecting is not settled.

\flushpar
We believe that the second case above never happens; conjugacy classes are
contained in analytic submanifolds. Unfortunately we were not able to
 prove this. Showing that the restriction of $\psi$ to $p\circ \pi([f])$
is a $C^1$ function
is equivalent to showing that every
 invariant measure is obtained by projecting
the measure from the Markov extension.

%% file: plref.tex
\tolerance=3000

\bigskip
\centerline{\bf References}
\bigskip

\parindent=40pt
\item{[BKNS]} H. Bruin, G. Keller,T. Nowicki, and S. van Strien, 
{\it Fibonacci maps, Part I. Absorbing Cantor Sets}, Preprint 1994.
\item{[BGG]} C. Brodiscou, Ch. Gillot, and G. Gillot, {\it Variations
 d'entropie
d'une famille de transformations de l'intervalle unit{\'e}, unimodales et
lin\'eaires par morceaux}, C.R.Acad.Sci. Paris {\bf 300} s{\'e}rie I (1985) 585--588.
\item{[L]}  M. Lyubich, {\it Combinatorics, geometry and attractors of
quasi-quadratic maps}, Preprint  SUNYSB \#18 (1992). 
\item{[LM]} M. Lyubich and J. Milnor, {\it The Fibonacci Unimodal map},
J.Amer.Math.Soc. {\bf 6} (1993) 425--457.
\item{[LY]} A. Lasota and Y. Yorke, {\it On the existence of
 invariant measures 
for piecewise monotone transformations}, Trans.Amer.Math.Soc. {\bf 183} (1973) 481--488. 
\item{[M]}  M. Martens, {\it Distortion results and invariant Cantor sets
of unimodal maps}, To appear in Erg.Th. \& Dyn.Sys.
\item{[MMS]} M. Martens, W. de Melo, and S.van Strien, {\it Julia-Fatou-Sullivan
theory for real one-dimensional dynamics}, Acta. Math. {\bf 168} (1992) 273--318.
\item{[Mi]} M. Misiurewicz, {\it Absolutely continuous invariant measures for 
certain maps of the interval}, Publ.Math. IHES {\bf 53} (1981) 17--51. 
\item{[MiV]} M. Misiurewicz and E. Visinescou, {\it Kneading sequences of skew
tent maps}, Ann.Inst. Henri Poincar\'e {\bf 27} (1991) 125--140.